\documentclass[a4paper,12pt]{article}
\usepackage{amsfonts,amssymb,latexsym,amsmath}
\usepackage{multicol}
\usepackage[cp1251]{inputenc}
\usepackage[russian]{babel}
\begin{document}
\newcounter{bnomer}
\newcounter{snomer}
\newcounter{diagram}
\setcounter{bnomer}{0} \setcounter{diagram}{0}
\renewcommand{\thesnomer}{\thebnomer.\arabic{snomer}}
\renewcommand{\thebnomer}{\arabic{bnomer}}

\newcommand{\sect}[1]{%
\setcounter{snomer}{0} \refstepcounter{bnomer}
\begin{center}\large{\textbf{\S \thebnomer.{ #1}}}\end{center}}

\newcommand{\thenv}[2]{%
\refstepcounter{snomer}
\par\addvspace{\medskipamount}\textbf{{#1} \thesnomer.}
{#2}\par\addvspace{\medskipamount}}

%====================================================

\renewcommand{\refname}{References}

\date{}
\title{The algebra of invariants of the adjoint action of the
unitriangular group in the nilradical of a~parabolic subalgebra}
\author{V. V. Sevostyanova\thanks{This research was partially
supported by the RFBR (project 08-01-00151-а)}}

\maketitle

\begin{center}
\parbox[b]{330pt}{\small\textsc{Abstract.} In the paper the
algebra of invariants of the adjoint action of the unitriangular
group in the nilradical of a parabolic subalgebra is studied. We set
up a conjecture on the structure of the algebra of invariants. The
conjecture is proved for parabolic subalgebras of special types.}
\end{center}

\vspace{0.5cm}

Consider the general linear group $\mathrm{GL}(n,K)$ defined over an
algebraically closed field $K$ of characteristic 0. Let $B$ ($N$,
respectively) be its Borel (maximal unipotent, respectively)
subgroup, which consists of triangular matrices with nonzero (unit,
respectively) elements on the diagonal. We fix a parabolic subgroup
$P$ that contains $B$. Denote by $\mathfrak{p}$, $\mathfrak{b}$ and
$\mathfrak{n}$ the Lie subalgebras in $\mathfrak{gl}(n,K)$ that
correspond to $P$, $B$ and $N$, respectively. We represent
$\mathfrak{p}=\mathfrak{r}\oplus\mathfrak{m}$ as the direct sum of
the nilradical $\mathfrak{m}$ and a block diagonal subalgebra
$\mathfrak{r}$ with sizes of blocks $(n_1,\ldots, n_s)$. The
subalgebra $\mathfrak{m}$ is invariant relative to the adjoint
action of the group $P$, therefore $\mathfrak{m}$ is invariant
relative to the adjoint action of the subgroups $B$ and $N$. We
extend this action to the representation in the algebra
$K[\mathfrak{m}]$ and in the field $K(\mathfrak{m})$. The subalgebra
$\mathfrak{m}$ contains a Zariski-open $P$-orbit, which is called
the \emph{Richardson orbit} (see~\cite{R}). Consequently,
$K[\mathfrak{m}]^P=K$. In this paper, we study the structure of the
algebra of invariants $K[\mathfrak{m}]^N$. In the case $P=B$, the
algebra of invariants $K[\mathfrak{m}]^N$ is the polynomial algebra
$K[x_{12}, x_{23},\ldots,x_{n-1,n}]$. Let $\mathfrak{r}$ be the sum
of two blocks; this case was considered by M. Brion in the
paper~\cite{B}. A complete description of the field of invariants
$K(\mathfrak{m})^N$ for any parabolic subalgebra is a result
of~\cite{S}. The question concerning the structure of the algebra of
invariantss $K[\mathfrak{m}]^N$ remains open and seems to be a
considerable challenge. We do not know when the algebra of
invariants $K(\mathfrak{m})^N$ is finitely generated.

In this paper, we consider the series of polynomials (see the
notations (\ref{A}), (\ref{B}), (\ref{C})) in the algebra of
invariants. We show that these polynomials generate the algebra of
invariants in special cases (Theorems \ref{Th_2k2}
and~\ref{Th1221}). We state a conjecture on the structure of the
algebra of invariants $K[\mathfrak{m}]^N$
(Conjecture~\ref{conjecture-main}).

%========================================================

\sect{The main definitions}

We begin with definitions. Every positive root $\gamma$ in
$\mathfrak{gl}(n,K)$ has the form~(see \cite{GG})
$\gamma=\varepsilon_i-\varepsilon_j$, $1\leqslant i<j\leqslant n$.
We identity a root $\gamma$ with the pair $(i,j)$ and the set of the
positive roots $\Delta^{\!+}$ with the set of pairs $(i,j)$, $i<j$.
The system of positive roots $\Delta^{\!+}_\mathfrak{r}$ of the
reductive subalgebra $\mathfrak{r}$ is a subsystem in
$\Delta^{\!+}$.

Let $\{E_{i,j}:~i<j\}$ be the standard basis in $\mathfrak{n}$. By
$E_\gamma$ denote the basis element $E_{i,j}$, where $\gamma=(i,j)$.

We define a relation in $\Delta^{\!+}$ such that
$\gamma'\succ\gamma$ whenever
$\gamma'-\gamma\in\Delta^{\!+}_\mathfrak{r}$. If
$\gamma\prec\gamma'$ or $\gamma\succ\gamma'$, then the roots
$\gamma$ and $\gamma'$ are said to be \emph{comparable}. Denote by
$M$ the set of $\gamma\in\Delta^{\!+}$ such that
$E_\gamma\in\mathfrak{m}$. We identify the algebra $K[\mathfrak{m}]$
with the polynomial algebra in the variables $x_{i,j}$, ~$(i,j)\in
M$.

\thenv{Definition}{A subset $S$ in $M$ is called a \emph{base} if
the elements in $S$ are not pairwise comparable and for any
$\gamma\in M\setminus S$ there exists $\xi\in S$ such that
$\gamma\succ\xi$.}

\thenv{Definition}{Let $A$ be a subset in $S$. We say that $\gamma$
is a \emph{minimal element} in $A$ if there is no $\xi\in A$ such
that $\gamma\succ\xi$.}

\medskip
Note that $M$ has a unique base $S$, which can be constructed in the
following way. We form the set $S_1$ of minimal elements in $M$. By
definition, $S_1\subset S$. Then we form a set $M_1$, which is
obtained from $M$ by deleting $S_1$ and all elements
$$\{\gamma\in M:\exists\ \xi\in S_1,\ \gamma\succ\xi\}.$$
The set of minimal elements $S_2$ in $M_1$ is also contained in $S$,
and so on. Continuing the process, we get the base $S$.

\thenv{Definition}{An ordered set of positive roots
$\{\gamma_1,\ldots,\gamma_s\}$ is called a \emph{chain} if
$\gamma_1=(a_1,a_2)$, $\gamma_2=(a_2,a_3)$, $\gamma_3=(a_3,a_4)$,
and so on. The number $s$ is called the \emph{length of a chain}.}

\thenv{Definition}{We say that two roots $\xi,\xi'\in S$ form an
\emph{admissible pair} $q=(\xi,\xi')$ if there exists
$\alpha_q\in\Delta^{\!+}_\mathfrak{r}$ such that the ordered set of
roots $\{\xi,\alpha_q,\xi'\}$ is a chain. Note that the root
$\alpha_q$ is uniquely determined by $q$.}

We form the set $Q:=Q(\mathfrak{p})$ that consists of admissible
pairs of roots in $S$. For every admissible pair $q=(\xi,\xi')$ we
construct a positive root $\varphi_q=\alpha_q+\xi'$. Consider the
subset $\Phi=\{\varphi_q:~ q\in Q\}$ in the set of positive roots.

\thenv{Definition}{The set $S\cup\Phi$ is called an \emph{expanded
base}.}

Using the given parabolic subgroup, we construct a diagram, which is
a square matrix in which the roots from $S$ are marked by the symbol
$\otimes$ and the roots from $\Phi$ are labeled by the symbol
$\times$. The other entries in the diagram are empty.

\thenv{Example}{Below a diagram for a parabolic subalgebra with
sizes of its diagonal blocks $(2,1,3,2)$ is given.}
\begin{center}\refstepcounter{diagram}
{\begin{tabular}{|p{0.1cm}|p{0.1cm}|p{0.1cm}|p{0.1cm}|p{0.1cm}|p{0.1cm}|p{0.1cm}|p{0.1cm}|c}
\multicolumn{2}{l}{{\small 1\quad 2}}&\multicolumn{2}{l}{{\small
3\quad 4}}&\multicolumn{2}{l}{{\small 5\quad 6}}&
\multicolumn{2}{l}{{\small 7\quad 8}}\\
\cline{1-8} \multicolumn{2}{|l|}{1}&&&$\otimes$&&&&{\small 1}\\
\cline{3-8} \multicolumn{2}{|r|}{1}&$\otimes$&&&&&&{\small 2}\\
\cline{1-8} \multicolumn{2}{|c|}{}&1&$\otimes$&&&&&{\small 3}\\
\cline{3-8} \multicolumn{3}{|c|}{}&\multicolumn{3}{|l|}{1}&$\times$&$\times$&{\small 4}\\
\cline{7-8} \multicolumn{3}{|c|}{}&\multicolumn{3}{|c|}{1}&$\times$&$\otimes$&{\small 5}\\
\cline{7-8} \multicolumn{3}{|c|}{}&\multicolumn{3}{|r|}{1}&$\otimes$&&{\small 6}\\
\cline{4-8} \multicolumn{6}{|c|}{}&\multicolumn{2}{|l|}{1}&{\small 7}\\
\multicolumn{6}{|c|}{}&\multicolumn{2}{|r|}{1}&{\small 8}\\
\cline{1-8} \multicolumn{8}{c}{Diagram \arabic{diagram}}\\
\end{tabular}}
\end{center}

Consider the formal matrix $\mathbb{X}$ in which the variables
$x_{i,j}$ occupy the positions $(i,j)\in M$ and the other entries
are equal to zero. For any root $\gamma=(a,b)\in M$ we denote by
$S_\gamma$ the set of $\xi=(i,j)\in S$ such that $i>a$ и $j<b$. Let
$S_\gamma=\{(i_1,j_1),\ldots,(i_k,j_k)\}$. Denote by $M_\gamma$ a
minor $M_I^J$ of the matrix $\mathbb{X}$ with ordered systems of
rows $I$ and columns $J$, where
$$I=\mathrm{ord}\{a,i_1,\ldots,i_k\},\quad J=\mathrm{ord}\{j_1,\ldots,j_k, b\}.$$

For every admissible pair $q=(\xi,\xi')$ such that $q$ corresponds
to $\varphi\in\Phi$, we construct the polynomial
\begin{equation}
L_\varphi=\sum_{\scriptstyle\alpha_1,\alpha_2\in\Delta^{\!+}_\mathfrak{r}\cup\{0\}
\atop\scriptstyle\alpha_1+\alpha_2=\alpha_q}
M_{\xi+\alpha_1}M_{\alpha_2+\xi'}.\label{L_q}
\end{equation}

\thenv{Theorem\label{M-L_independ}}{\cite{PS} \emph{For an arbitrary
parabolic subalgebra}, \emph{the system of polynomials}
$$\{M_\xi,~\xi\in S,~L_{\varphi},~\varphi\in\Phi,\}$$
\emph{is contained in $K[\mathfrak{m}]^N$ and is algebraically
independent over $K$}.}

Denote by $\mathcal{Y}$ the subset in $\mathfrak{m}$ that consists
of matrices of the form
$$\sum_{\xi\in S}c_{\xi}E_\xi+
\sum_{\varphi\in\Phi}c'_{\varphi}E_\varphi,$$ where $c_{\xi}\neq0$
and $c'_{\varphi}\neq0$.

\thenv{Definition}{The matrices from $\mathcal{Y}$ are said to be
\emph{canonical}.}

The proof of the following theorem is found in~\cite{S}.

\thenv{Theorem\label{Exist_of_representative}}{\emph{There exists a
nonempty Zariski-open subset $U\subset\mathfrak{m}$ such that the
$N$-orbit of any $x\in U$ intersects $\mathcal{Y}$ at a unique
point}.}

Now let $\mathcal{S}$ be the set of denominators generated by minors
$M_\xi$, $\xi\in S$. Consider the localization
$K[\mathfrak{m}]^N_\mathcal{S}$ of the algebra of invariants
$K[\mathfrak{m}]^N$ with respect to $\mathcal{S}$. Since the minors
$M_\xi$ are $N$-invariants, we have
$$K[\mathfrak{m}]^N_\mathcal{S}=\big(K[\mathfrak{m}]_S\big)^N.$$

The following results are consequences of
Theorem~\ref{Exist_of_representative}.

\thenv{Theorem\label{Th_local_field}}{\emph{The ring
$K[\mathfrak{m}]^N_S$ is the ring of polynomials in} $M_\xi^{\pm
1}$, $\xi\in S$, \emph{and in} $L_{\varphi}$, $\varphi\in\Phi$.}

\thenv{Theorem\label{Th_invariant_field}}{\emph{The field of
invariants $K(\mathfrak{m})^N$ is the field of rational functions
of} $M_\xi$, $\xi\in S$, \emph{and} $L_{\varphi}$,
$\varphi\in\Phi$.}

The polynomials $M_{\xi}$ and $L_{\varphi}$ do not generate the
algebra of invariants. We prove in \S4 and \S5, respectively, that
if the reductive subalgebra~$\mathfrak{r}$ consists of three blocks
or of four blocks with sizes $(1,2,2,1)$, then the algebra of
invariants $K[\mathfrak{m}]^N$ is not generated by the invariants
$M_{\xi}$ and $L_{\varphi}$. Let the reductive subalgebra consist of
blocks with sizes $(2,k,2)$, $k>3$, or with sizes $(1,2,2,1)$. We
construct a polynomial $D$ such that the algebra of invariants is
generated by $D$ and by the invariants $M_{\xi}$ and $L_{\varphi}$.

%============================================================
\sect{The other definition of the invariant $L_{\varphi}$}

In the present chapter, we give the alternative definition of the
invariant~$L_{\varphi}$, where $\varphi\in\Phi$.

Let $I$ and $I'$ be two systems of rows, and let $J$ and $J'$ be two
systems of columns:
$$I=\{i,i+1,\ldots,i+k\},$$
$$J=\{j,j+1,\ldots,j+l\},$$
$$I'=\{i',i'+1,\ldots,i'+k'\},$$
$$J'=\{j',j'+1,\ldots,j'+l'\},$$
where $i+k<j'\leqslant i'$ and $j'+l'\leqslant i'+k'<j$. Suppose the
minors $\mathbb{X}_{I}^{J'}$ and $\mathbb{X}_{I'}^{J}$ are bordered
by zeros from below and left in~$\mathbb{X}$, where
$\mathbb{X}_{I}^{J'}$ is the minor of the matrix~$\mathbb{X}$ with
systems of rows $I$ and columns $J'$. In other words, if for any
positive root $(a,b)$ one of the following conditions holds
\begin{itemize}
\item[1)] $a\in I$ and $b<j'$,
\item[2)] $a\in I'$ and $b<j$,
\item[3)] $a>i+k$ and $b\in J'$,
\item[4)] $a>i'+k'$ and $b\in J$,
\end{itemize}
then the root $(a,b)$ is not in the system $M$. Let $k+k'=l+l'$. The
following determinant is one of order $k+k'+2$:
\begin{equation}
D_{I,I'}^{J,J'}=\left|\begin{tabular}{c|c}
\raisebox{5pt}[25pt]{$\mathbb{X}_I^{J'}$}&
\raisebox{5pt}[25pt]{$(\mathbb{X}^2)_I^J$}\\
\hline \raisebox{0pt}[22pt]{\Large$0$}&
\raisebox{0pt}[22pt]{$\mathbb{X}_{I'}^J$}\\
\end{tabular}
\right|,\label{Sbornii-minor}
\end{equation}
where the minor $\mathbb{X}_I^J$ of the matrix $\mathbb{X}$ is
formed by systems of rows~$I$ and columns~$J$.

\thenv{Definition}{The determinant $D_{I,I'}^{J,J'}$ is called a
\emph{combined minor}.}

\thenv{Lemma\label{Minor-Inv}}{\emph{A combined
minor}~(\ref{Sbornii-minor}) \emph{is an $N$-invariant}.}

\textsc{Proof.} A polynomial $H$ is an $N$-invariant if $H$ is an
invariant under the adjoint action of all one-parameter subgroups
$g_m(t)=E+tE_{m,m+1}$, $1\leqslant m<n$, where $E$ is a matrix unit.
The adjoint action on the matrix $\mathbb{X}$ by $g_m(t)$ reduces to
the composition of two transformations:
\begin{itemize}
\item[1)] the row with number $m+1$ multiplied by~$t$ is added to the
row $m$ of the matrix $\mathbb{X}$;
\item[2)]  the column with the number $m$ multiplied by $-t$ is
added to the column with the number $m+1$ of the matrix
$\mathbb{X}$.
\end{itemize}

We show that the polynomial~(\ref{Sbornii-minor}) is an invariant
under the adjoint action of $g_{m,m+1}(t)$. Without loss of
generality, it can be assumed that $i+k=j'-1$, $i'+k'=j-1$ and the
minors $\mathbb{X}_I^{J'}$, $\mathbb{X}_{I'}^J$, and
$(\mathbb{X}^2)_I^J$ have the following form
$$\mathbb{X}_I^{J'}=\left|\begin{array}{cccc}
x_{i,j'}&x_{i,j'+1}&\ldots&x_{i,j'+l'}\\
x_{i+1,j'}&x_{i+1,j'+1}&\ldots&x_{i+1,j'+l'}\\
\ldots&\ldots&\ldots&\ldots\\
x_{i+k,j'}&x_{i+k,j'+1}&\ldots&x_{i+k,j'+l'}
\end{array}\right|,$$
$$\mathbb{X}_{I'}^{J}=\left|\begin{array}{cccc}
x_{i',j}&x_{i',j+1}&\ldots&x_{i',j+l}\\
x_{i'+1,j}&x_{i'+1,j+1}&\ldots&x_{i'+1,j+l}\\
\ldots&\ldots&\ldots&\ldots\\
x_{i'+k',j}&x_{i'+k',j+1}&\ldots&x_{i'+k',j+l}
\end{array}\right|,$$
$$(\mathbb{X}^2)_I^J=\left|\begin{array}{cccc}
\widetilde{x}_{i,j}&\widetilde{x}_{i,j+1}&\ldots&\widetilde{x}_{i,j+l}\\
\widetilde{x}_{i+1,j}&\widetilde{x}_{i+1,j+1}&\ldots&\widetilde{x}_{i+1,j+l}\\
\ldots&\ldots&\ldots&\ldots\\
\widetilde{x}_{i+k,j}&\widetilde{x}_{i+k,j+1}&\ldots&\widetilde{x}_{i+k,j+l}
\end{array}\right|,$$
where
$\displaystyle\widetilde{x}_{a,b}=\sum_{p=j'}^{i'+k'}x_{a,p}x_{p,b}$
for any $a\in I$ и $b\in J$. If $m<i$ or $m+1>j+l$, then the adjoint
action of the element $g_{m,m+1}(t)$ doesn't change the
minor~(\ref{Sbornii-minor}).

\begin{enumerate}
\item
Suppose $i\leqslant m<i+k$; then the action of the element
$g_{m,m+1}(t)$ changes the variables $x_{a,b}$ of the minor
$D_{I,I'}^{J,J'}$ such that $a=m$ and doesn't change the other
variables. Therefore
$$g_{m,m+1}(t)\,x_{m,p}=x_{m,p}-tx_{m+1,p},\ j'\leqslant p\leqslant i'+k';$$
$$g_{m,m+1}(t)\sum_{p=j'}^{i'+k'}x_{m,p}x_{p,b}=
\sum_{p=j'}^{i'+k'}x_{m,p}x_{p,b}-t\sum_{p=j'}^{i'+k'}x_{m+1,p}x_{p,b}.$$
Then we have
$$g_{m,m+1}(t)\mathbb{X}_{I'}^J=\mathbb{X}_{I'}^J,$$
$$g_{m,m+1}(t)\mathbb{X}_I^{J'}=
\mathbb{X}_I^{J'}-tD_1.$$ The order of $D_1$ is equal to the order
of $\mathbb{X}_{I}^{J'}$. The minor $D_1$ is formed by the
columns~$J'$ and $D_1$ has two identical rows with numbers~$m+1$
and~$m$. Further,
$$g_{m,m+1}(t)(\mathbb{X}^2)_I^{J}=
(\mathbb{X}^2)_I^{J}-tD_2,$$ where the order of $D_2$ is equal to
the order of $(\mathbb{X}^2)_{I}^{J}$. The minor $D_2$ has two
identical rows with numbers~$m$ and~$m+1$.

Thus,
$$g_{m,m+1}(t)D_{I,I'}^{J,J'}=g_{m,m+1}(t)\left|\begin{tabular}{c|c}
\raisebox{5pt}[25pt]{$\mathbb{X}_I^{J'}$}&
\raisebox{5pt}[25pt]{$(\mathbb{X}^2)_I^J$}\\
\hline \raisebox{0pt}[22pt]{\Large$0$}&
\raisebox{0pt}[22pt]{$\mathbb{X}_{I'}^J$}\\
\end{tabular}
\right|=$$$$=\left|\begin{tabular}{c|c}
\raisebox{5pt}[25pt]{$\mathbb{X}_I^{J'}$}&
\raisebox{5pt}[25pt]{$(\mathbb{X}^2)_I^J$}\\
\hline \raisebox{0pt}[22pt]{\Large$0$}&
\raisebox{0pt}[22pt]{$\mathbb{X}_{I'}^J$}\\
\end{tabular}
\right|-t\left|\begin{tabular}{c|c} \raisebox{5pt}[25pt]{$D_1$}&
\raisebox{5pt}[25pt]{$D_2$}\\
\hline \raisebox{0pt}[22pt]{\Large$0$}&
\raisebox{0pt}[22pt]{$\mathbb{X}_{I'}^J$}\\
\end{tabular}
\right|.$$ The minor that consists of $D_1$, $D_2$, and
$\mathbb{X}_{I'}^J$ has two identical rows; therefore this minor is
equal to zero. We have
$$g_{m,m+1}(t)D_{I,I'}^{J,J'}=D_{I,I'}^{J,J'}.$$

\item
Let $i+k\leqslant m<j$. For any $a\in I$ and $b\in J$ we have
$$g_{m,m+1}(t)\sum_{p=j'}^{i'+k'}x_{a,p}x_{p,b}=
\sum_{p=j'}^{m-1}x_{a,p}x_{p,b}+x_{a,m}(x_{m,b}-tx_{m+1,b})+$$$$
+(x_{a,m+1}+tx_{a,m})x_{m+1,b}+\sum_{p=m+2}^{i'+k'}x_{a,p}x_{p,b}
=\sum_{p=j'}^{i'+k'}x_{a,p}x_{p,b},$$ i.e., the adjoint action of
the element $g_{m,m+1}(t)$ does not change the minor
$(\mathbb{X}^2)_I^J$. Besides, the action of the element
$g_{m,m+1}(t)$ adds the column with the number $m$ to the $(m+1)$th
column of the minor~$\mathbb{X}_I^{J'}$ and the $(m+1)$th row to the
$m$th row of the minor~$\mathbb{X}_{I'}^J$. Thus the minor
$D_{I,I'}^{J,J'}$ is a $g_{m,m+1}(t)$-invariant.

\item For $j\leqslant m<j+l$, the proof is similarly.~$\Box$
\end{enumerate}

We shell show that the invariant~$L_{\varphi}$ is a combined
minor~(\ref{Sbornii-minor}) for any systems~$I,I',J,J'$.

\thenv{Definition}{Any root $(i,j)\in M$ satisfying the conditions
$$R_{k-1}<i\leqslant R_k\mbox{ and }R_k<j\leqslant n,$$
is called a \emph{root lying to the right of the $k$th block} in
$\mathfrak{r}$.}

\thenv{Definition}{Any root $(i,j)\in M$ satisfying the conditions
$$1\leqslant i\leqslant R_{k-1}\mbox{ and }R_{k-1}<j\leqslant R_{k},$$
is called a \emph{root lying to the above of the $k$th block} in
$\mathfrak{r}$.}

Let $\varphi$ be any root from $\Phi$. The root~$\varphi$
corresponds to an admissible pair $(\xi,\xi')$. Let $\xi=(i,j)$,
$\xi'=(l,m)$, and $i<l$. Assume that $\varphi$ lies to the right of
the $k$th block in $\mathfrak{r}$, i.e., $R_{k-1}<j\leqslant R_k$.
Let the set $S_{\xi}$ has $p$ roots and let the set $S_{\xi'}$ has
$q$ roots. We have $i+p=R_{k-1}$, i.e., the minimal root $\gamma$ in
$S_{\xi}$, which has the greatest number of row, lies in the row
with number $R_{k-1}$. Similarly, the root of the base that has the
least number of  column is also minimal; so we have $l-q=R_k+1$.

Denote
$$I=\{i,i+1,\ldots,i+p\},\ J=\{m-q,m-q+1,\ldots,m\},$$
$$I'=\{l+1,l+2,\ldots,l+q\},\ J'=\{j-p,j-p+1\ldots,j-1\}.$$
Then $|I|=p+1,\ |J|=q+1,\ |I'|=q,\ |J'|=p.$ Consider the following
combined minor
$$\widetilde{L}_{\varphi}=D_{I,I'}^{J,J'}=\left|\begin{tabular}{c|c}
\raisebox{5pt}[25pt]{$\mathbb{X}_I^{J'}$}&
\raisebox{5pt}[25pt]{$(\mathbb{X}^2)_I^J$}\\
\hline \raisebox{0pt}[22pt]{\Large$0$}&
\raisebox{0pt}[22pt]{$\mathbb{X}_{I'}^J$}\\
\end{tabular}
\right|,$$ where $\mathbb{X}_I^J$ is the minor of $\mathbb{X}$ with
the systems of rows~$I$ and columns~$J$. The order of
$\widetilde{L}_{\varphi}$ is equal to $(p+q+1)$.

\thenv{Proposition\label{Def_L_phi}}{\emph{We have the equation}
$\widetilde{L}_{\varphi}=L_{\varphi}$.}

\textsc{Proof.} By Lemma~\ref{Minor-Inv}, the
polynomial~$\widetilde{L}_{\varphi}$ is an invariant. Let
$$p:\,K[\mathfrak{m}]^N\rightarrow K[\mathcal{Y}],\quad
p(f)=f|_{\mathcal{Y}}$$ be a restriction homomorphism. We calculate
the invariants~$L_{\varphi}$ and~$\widetilde{L}_{\varphi}$
on~$\mathcal{Y}$. Let $\mathbb{Y}$ be the formal matrix in which the
variables $x_{i,j}$ occupy the positions $(i,j)\in S\cup\Phi$ and
the other entries are equal to zero.

Assume that the root $\varphi$ corresponds to the admissible pair
$(\xi,\xi')$, where
$$\xi=(i,j),\ \xi'=(l,m),$$ and the sets $S_{\xi}$, $S_{\xi'}$
consist of the roots
$$S_{\xi}=\{\xi_1,\xi_2,\ldots,\xi_p\},\
S_{\xi'}=\{\xi'_1,\xi'_2,\ldots,\xi'_q\},$$ respectively.

Then using elementary transformations of columns, the minor
$\mathbb{Y}_I^{J'}$ of the matrix~$\mathbb{Y}$ is reduced to the
following minor of order $(p+1)\times p$:
$$\left|\begin{array}{cccc}
0&0&\ldots&0\\
0&0&\ldots&x_{\xi_p}\\
\ldots&\ldots&\ldots&\ldots\\
0&x_{\xi_2}&\ldots&0\\
x_{\xi_1}&0&\ldots&0\\
\end{array}\right|.$$
The minor $\mathbb{Y}_{I'}^{J}$ is reduced to the following minor of
order $q\times(q+1)$ by elementary transformations of rows:
$$\left|\begin{array}{ccccc}
0&0&\ldots&x_{\xi'_q}&0\\
\ldots&\ldots&\ldots&\ldots&\ldots\\
0&x_{\xi'_2}&\ldots&0&0\\
x_{\xi'_1}&0&\ldots&0&0\\
\end{array}\right|.$$

We have
$$\widetilde{L}_{\varphi}=\pm
\left|\begin{array}{ccccccccc}
0&0&\ldots&0&&&&&\\
0&0&\ldots&x_{\xi_p}&&&&\\
\ldots&\ldots&\ldots&\ldots&&\multicolumn{3}{c}{(\mathbb{Y}^2)_I^J}&\\
0&x_{\xi_2}&\ldots&0&&&&&\\
x_{\xi_1}&0&\ldots&0&&&&&\\
&&&&\makebox[0pt][r]{0\ }&0&\ldots&x_{\xi'_q}&0\\
&&&&\makebox[0pt][r]{$\ldots$}&\ldots&\ldots&\ldots&\ldots\\
&0&&&\makebox[0pt][r]{0\ }&x_{\xi'_2}&\ldots&0&0\\
&&&&\makebox[0pt][r]{$x_{\xi'_1}$}&0&\ldots&0&0\\
\end{array}\right|=$$
\begin{equation}
=\pm
x_{\xi}x_{\varphi}\prod_{s=1}^{p}x_{\xi_s}\cdot\prod_{t=1}^{q}x_{\xi'_t}
.\label{minor_L}
\end{equation}

Direct calculations show that the invariant $L_{\varphi}$
on~$\mathbb{Y}$ is equal to~(\ref{minor_L}).

Let we show that the mapping~$p$ is an embedding. Indeed, if
$f\in\mathrm{Ker}\,p$, then
$f\left(\mathrm{Ad}_N\,\mathcal{Y}\right)=0$. By
Theorem~\ref{Exist_of_representative}, we have
$\mathrm{Ad}_N\,\mathcal{Y}$. Then $f\equiv0$ and $p$ is an
embedding.

The invariant $L_{\varphi}$ is equal to $\widetilde{L}_{\varphi}$
on~$\mathcal{Y}$. Therefore $L_{\varphi}$ is equal to
$\widetilde{L}_{\varphi}$ on the pre-image $p^{-1}(\mathcal{Y})$.
Finally, we get $L_{\varphi}=\widetilde{L}_{\varphi}$.~$\Box$

%========================================================

\sect{The additional set of $N$-invariants}

Let $\mathfrak{p}$ be any parabolic subalgebra.

Now we introduce a system of $N$-invariants that are not in the
algebra $K[M_{\xi}, L_{\varphi}]_{{\xi\in S}\atop {\varphi\in
\Phi}}$.

\medskip
Let $(i,j)\in S\cup\Phi$. Denote
$$\begin{array}{r}L_{i,j}=\left\{\begin{array}{l}
L_{(i,j)}\mbox{ if }(i,j)\in\Phi;\\
M_{(a,i)}\cdot M_{(i,j)}\mbox{ if }(i,j)\in S\mbox{ and there is a
number }a\mbox{ such that }\\
\end{array}\right.\\
\hspace{9.5cm}(a,i)\in S.
\end{array}$$

\medskip
Assume that for any numbers $m<l$ and $i<j$ the roots
\begin{equation}
(m,i),(l,i),(m,j),(l,j)\label{4roots}
\end{equation}
are contained in the expanded base. Obviously, the root $(m,i)$ lies
in the system~$\Phi$. Suppose that the root $(m,i)$ corresponds to
some admissible pair $(\xi,\xi')\in Q$, $\xi,\xi'\in S$.

\thenv{Notation}{Denote
\begin{equation}
A_{m}^{i,j}=\frac{L_{m+1,i}L_{m,j}-L_{m+1,j}L_{m,i}}{M_{\xi}}.\label{A}
\end{equation}}

Now assume that the following condition holds for the
roots~(\ref{4roots}): there is not a number $\widetilde{i}$,
$i<\widetilde{i}<j$, such that $(m,\widetilde{i})\in S\cup\Phi$.

\thenv{Notation}{Denote
\begin{equation}
B^{\,i}_{m,l}=\frac{L_{m,j}L_{l,i}-L_{l,j}L_{m,i}}{M_{\xi'}},\label{B}
\end{equation}
\begin{equation}
C_{m}^{\,i}=\frac{L_{m+1,i}L_{m,j}-L_{m+1,j}L_{m,i}}{M_{\xi}\cdot
M_{\xi'}}.\label{C}
\end{equation}}

Clearly, all rational function (\ref{A}), (\ref{B}), (\ref{C}) are
contained in the localization $K[\mathfrak{m}]^N_{\mathcal{S}}$,
where the set $\mathcal{S}$ is generated by the minors $M_{\xi}$,
$\xi\in S$. We shall prove that the functions (\ref{A}), (\ref{B}),
(\ref{C}) are polynomials.

Let the root $(m,i)$ be to the right of the $k$th block in
$\mathfrak{r}$, i.e., $R_{k-1}<m\leqslant R_k$. Let the
invariant~$C_m^{\,i}$ is constructed by roots
$$(m,i),(m+1,i),(m,j),(m+1,j),$$
which contained in the expanded base. Then there exists no
number~$\widetilde{i}$ such that $i<\widetilde{i}<j$ and
$(m,\widetilde{i})\in\Phi$. If the root~$(m+1,j)$ lies in the
system~$\Phi$, then $(m+1,j)$ corresponds to some admissible pair
$(\xi,\xi')$. If $(m+1,j)\in S$, then there exists a root $\xi$ in
the column with number $m+1$. Indeed, $(m+1,i)\in\Phi$ corresponds
to an admissible pair that contains $\xi$. In any case, we have that
there exist the roots $\xi=(a,m+1)$ and $\xi'=(b,j)$ for some
numbers~$a$ and~$b$. Consider the sets $S_{\xi}$ and $S_{\xi'}$.
Suppose $S_{\xi}$ contains $p$ roots and $S_{\xi'}$ contains $q$
roots. Since the root~$\xi$ lies to the above of the $k$th block
in~$\mathfrak{r}$, then we have $a+p=R_{k-1}$, i.e., the root
from~$S_{\xi}$ that has the greatest number of row is contained in
the row~$R_{k-1}$. Similarly, $j-q=R_k+1$.

\thenv{Example\label{314123}}{Suppose the reductive
subalgebra~$\mathfrak{r}$ is formed by the blocks $(3,1,4,1,2,3).$
We have the following diagram.
\begin{center}\refstepcounter{diagram}
{\begin{tabular}{|p{0.1cm}|p{0.1cm}|p{0.1cm}|p{0.1cm}|p{0.1cm}|
p{0.1cm}|p{0.1cm}|p{0.1cm}|p{0.1cm}|p{0.1cm}|p{0.1cm}|p{0.1cm}|
p{0.1cm}|p{0.1cm}|c} \multicolumn{2}{l}{{\small 1\hspace{5pt}
2}}&\multicolumn{2}{l}{{\small 3\hspace{5pt}
4}}&\multicolumn{2}{l}{{\small 5\hspace{5pt} 6}}&
\multicolumn{2}{l}{{\small 7\hspace{5pt}
8}}&\multicolumn{2}{c}{{\small\makebox[0.2cm][c]{9\hspace{5pt}
10}}}&\multicolumn{2}{c}{{\small\makebox[0.2cm][c]{11\hspace{4pt}12}}}&
\multicolumn{2}{c}{{\small\makebox[0.2cm][c]{13\hspace{4pt}14}}}\\
\cline{1-14} \multicolumn{3}{|l|}{1}&&&&$\otimes$&&&&&&&&{\small 1}\\
\cline{4-14} \multicolumn{3}{|c|}{1}&&&$\otimes$&&&&&&&&&{\small 2}\\
\cline{4-14} \multicolumn{3}{|r|}{\hspace{30pt}1}&$\otimes$&&&&&&&&&&&{\small 3}\\
\cline{1-14} \multicolumn{3}{|l|}{}&1&$\otimes$&&&&&&&&&&{\small 4}\\
\cline{4-14} \multicolumn{4}{|c|}{}&\multicolumn{4}{|l|}{1}&
$\times$&&$\times$&&&$\times$&{\small 5}\\
\cline{9-14}
\multicolumn{4}{|c|}{}&\multicolumn{4}{|l|}{\hspace{15pt}1}&
$\times$&&$\times$&&&$\otimes$&{\small 6}\\
\cline{9-14}
\multicolumn{4}{|c|}{}&\multicolumn{4}{|l|}{\hspace{30pt}1}&
$\times$&&$\otimes$&&&&{\small 7}\\
\cline{9-14}
\multicolumn{4}{|c|}{}&\multicolumn{4}{|l|}{\hspace{45pt}1}&
$\otimes$&&&&&&{\small 8}\\
\cline{5-14} \multicolumn{8}{|c|}{}&1&$\otimes$&&&&&{\small 9}\\
\cline{9-14}
\multicolumn{9}{|c|}{}&\multicolumn{2}{|l|}{1}&$\times$&$\otimes$&&{\small 10}\\
\cline{12-14}
\multicolumn{9}{|c|}{}&\multicolumn{2}{|r|}{1}&$\otimes$&&&{\small 11}\\
\cline{10-14} \multicolumn{11}{|c|}{}&\multicolumn{3}{|l|}{1}&{\small 12}\\
\multicolumn{11}{|c|}{}&\multicolumn{3}{|c|}{1}&{\small 13}\\
\multicolumn{11}{|c|}{}&\multicolumn{3}{|r|}{1}&{\small 14}\\
\cline{1-14} \multicolumn{14}{c}{Diagram \arabic{diagram}}\\
\end{tabular}\label{maindiagram}}
\end{center}}
Consider the invariant~$C_{5}^{11}$. The invariant corresponds to
the roots $$(5,11),(6,11),(5,14),(6,14)$$ from the expanded base. By
the above notations, we have $\xi=(2,6)$ and $\xi'=(6,14)$, $p=2$,
and $q=5$.

Denote
$$I=\{a,a+1,\ldots,a+p\},\ J=\{j-q,j-q+1,\ldots,j\},$$
$$I'=\{b+2,b+3,\ldots,b+q\},\ J'=\{m-p+1,m-p+2\ldots,m-1\}.$$
Then $|I|=p+1,\ |J|=q+1,\ |I'|=q-1,\ |J'|=p-1.$ Consider the
following combined minor of order $(p+q)$:
$$\widetilde{C}_m^{\,i}=D_{I,I'}^{J,J'}=\left|\begin{tabular}{c|c}
\raisebox{5pt}[25pt]{$\mathbb{X}_I^{J'}$}&
\raisebox{5pt}[25pt]{$(\mathbb{X}^2)_I^J$}\\
\hline \raisebox{0pt}[22pt]{\Large$0$}&
\raisebox{0pt}[22pt]{$\mathbb{X}_{I'}^J$}\\
\end{tabular}
\right|,$$ where the minor $\mathbb{X}_I^J$ of the formal matrix is
formed by the intersection of rows~$I$ and columns~$J$.

\thenv{Example}{The polynomial~$\widetilde{C}_5^{11}$ of
Example~\ref{314123} is formed in the following way.
$$\left|\begin{tabular}{c|c}
\raisebox{5pt}[25pt]{$\mathbb{X}_{2,3,4}^{4}$}&
\raisebox{5pt}[25pt]{$(\mathbb{X}^2)_{2,3,4}^{9,10,11,12,13,14}$}\\
\hline \raisebox{0pt}[22pt]{\Large$0$}&
\raisebox{0pt}[22pt]{$\mathbb{X}_{8,9,10,11}^{9,10,11,12,13,14}$}\\
\end{tabular}
\right|.$$ The reader will easily prove that
$\widetilde{C}_5^{11}=C_5^{11}$.}

\thenv{Proposition\label{C=C_as-minor}}{\emph{We have
$C_m^{\,i}=\widetilde{C}_m^{\,i}$ for any numbers $m$ and}~$i$.}

\textsc{Proof.} Let us show that $\widetilde{C}_m^{\,i}=C_m^{\,i}$.

Let $$p:\,K[\mathfrak{m}]^N\rightarrow K[\mathcal{Y}],\quad
p(f)=f|_{\mathcal{Y}}$$ be the  the restriction homomorphism. We
show that $p$ is an embedding. Indeed, if $f\in\mathrm{Ker}\,p$,
then $f\left(\mathrm{Ad}_N\,\mathcal{Y}\right)=0$. Since by
Theorem~\ref{Exist_of_representative}, $\mathrm{Ad}_N\,\mathcal{Y}$
is contained a Zariski-open subset, then $f\equiv0$. Therefore $p$
is an embedding.

By Lemma~\ref{Minor-Inv}, the minor $\widetilde{C}_m^{\,i}$ is an
invariant. We show that the invariant~$L_{\varphi}$ is equal
to~$\widetilde{L}_{\varphi}$ on~$\mathcal{Y}$. Let $\mathbb{Y}$ be
the formal matrix such that the entries $(i,j)$ in $\mathbb{Y}$
contain the variables $x_{i,j}$ if $(i,j)\in S\cup\Phi$ and the
other entries is filled by zeros. We calculate the
polynomial~$\widetilde{C}_m^{\,i}$ and the invariant~$C_m^{\,i}$ on
$\mathbb{Y}$.

Now assume that the roots
$$\varphi_1=(m,i),\ \varphi_2=(m+1,i),\ \varphi_3=(m,j),\ \varphi_4=(m+1,j)$$
are in the expanded base and there exists no number $\widetilde{i}$
such that $i<\widetilde{i}<j$ and $(m,\widetilde{i})\in\Phi$.
Suppose the root $\varphi=(m+1,j)$ corresponds to the admissible
pair
$$\xi=(a,m+1),\ \xi'=(b,j).$$
Assume that the sets $S_{\xi}$ and $S_{\xi'}$ consist of the roots
$$S_{\xi}=\{\xi_1,\xi_2,\ldots,\xi_p\},\
S_{\xi'}=\{\xi'_1,\xi'_2,\ldots,\xi'_q\},$$ where the roots $\xi_p$
and $\xi'_q$ lie in the $m$th column and in the $(b+1)$th row,
respectively. Then by elementary transformations of columns, the
minor $\mathbb{Y}_I^{J'}$ of the matrix~$\mathbb{Y}$ is reduced to
the following minor of order $(p+1)\times(p-1)$:
$$\left|\begin{array}{cccc}
0&0&\ldots&0\\
0&0&\ldots&0\\
0&0&\ldots&x_{\xi_{p-1}}\\
\ldots&\ldots&\ldots&\ldots\\
0&x_{\xi_2}&\ldots&0\\
x_{\xi_1}&0&\ldots&0\\
\end{array}\right|,$$
By elementary transformations of rows, the minor
$\mathbb{Y}_{I'}^{J}$ is reduced to the following minor of order
$(q-1)\times(q+1)$:
$$\left|\begin{array}{cccccc}
0&0&\ldots&x_{\xi'_{q-1}}&0&0\\
\ldots&\ldots&\ldots&\ldots&\ldots&\ldots\\
0&x_{\xi'_2}&\ldots&0&0&0\\
x_{\xi'_1}&0&\ldots&0&0&0\\
\end{array}\right|.$$
We have
$$\widetilde{C}_{m}^{\,i}=\pm
\left|\begin{array}{cccccccccc}
0&0&\ldots&0&&&&&&\\
0&0&\ldots&0&&&&&&\\
0&0&\ldots&x_{\xi_{p-1}}&&\multicolumn{3}{c}{(\mathbb{Y}^2)_I^J}&&\\
\ldots&\ldots&\ldots&\ldots&&&&&\\
0&x_{\xi_2}&\ldots&0&&&&&\\
x_{\xi_1}&0&\ldots&0&&&&&\\
&&&&\makebox[0pt][r]{0\ }&0&\ldots&x_{\xi'_{q-1}}&0&0\\
&&&&\makebox[0pt][r]{$\ldots$}&\ldots&\ldots&\ldots&\ldots&\ldots\\
&0&&&\makebox[0pt][r]{0\ }&x_{\xi'_2}&\ldots&0&0&0\\
&&&&\makebox[0pt][r]{$x_{\xi'_1}$}&0&\ldots&0&0&0\\
\end{array}\right|=$$
$$=\pm\prod_{s=1}^{p-1}x_{\xi_s}\cdot\prod_{t=1}^{q-1}x_{\xi'_t}
\cdot\left|\begin{array}{cc}
x_{\xi}x_{\varphi_2}&x_{\xi}x_{\varphi_4}\\
x_{\xi_p}x_{\varphi_1}&x_{\xi_p}x_{\varphi_3}
\end{array}\right|=$$
\begin{equation}
=\pm x_{\xi}\prod_{s=1}^{p}x_{\xi_s}\cdot\prod_{t=1}^{q-1}x_{\xi'_t}
\cdot(x_{\varphi_2}x_{\varphi_3}-x_{\varphi_1}x_{\varphi_4}).\label{minor}
\end{equation}
By direct calculations, $C_m^{\,i}$ coincides with~(\ref{minor}) on
$\mathbb{Y}$. Thus $C_m^{\,i}$ coincides with~(\ref{minor})
on~$\mathcal{Y}$, and we have $\widetilde{C}_m^{\,i}$ is equal to
$C_m^{\,i}$ on the pre-image $p^{-1}(\mathcal{Y})$. So,
$C_m^{\,i}=\widetilde{C}_m^{\,i}$.~$\Box$

\medskip
Thus the invariant $C_m^{\,i}$ is a polynomial.

\thenv{Proposition\label{Pr_AB}}{\emph{We have} $A_m^{i,j},\
B^{\,i}_{m,l}\in K[\mathfrak{m}]^N.$}

\textsc{Proof.}Let us prove that $A_m^{i,j}$ is a polynomial. The
proof for $B^{\,i}_{m,l}$ is similarly.

Suppose the roots $(m,i)$, $(m,j)$, $(m+1,i)$, and $(m+1,j)$, $i<j$,
are contained in the expanded base and these roots lie to the right
of the $k$th block in the subalgebra $\mathfrak{r}$. We show that
$A_m^{i,j}$ is a polynomial. The proof is by induction.

Assume that there are $p$ roots to the above of the $k$th block
in~$S$. It is easily shown that these roots are the following ones:
$$(i_1,R_{k-1}+1),(i_2,R_{k-1}+2),\ldots,(i_p,R_{k-1}+p)$$
for some $i_1>i_2>\ldots>i_p$. Note that we have $i_1=R_{k-1}$ for
the root $(i_1,R_{k-1}+1)$. Suppose there are $q$ roots in~$S$ to
right of the $k$th block:
$$(R_k,j_1),(R_k-1,j_2),\ldots,(R_k-q+1,j_q)$$
for some $j_1<j_2<\ldots<j_q$. We have $j_1=R_k+1$. There exist the
numbers~$v$ and~$w$ such that $i=j_v$ and $j=j_w$. Assume that
$w-v>2$. Suppose the root $(m,j_w)$ corresponds to the admissible
pair~$(\xi,\xi')$. We fix the number $m$. We have the following
equation
\begin{equation}
L_{m+1,j_v}M_{\xi'}C_m^{j_{w-1}}=
A_m^{j_v,j_{w}}L_{m+1,j_{w-1}}-A_m^{j_v,j_{w-1}}L_{m+1,j_{w}}\label{equation}
\end{equation}
for any $1<w\leqslant q$. By Proposition~\ref{C=C_as-minor}, the
left part of the equation~(\ref{equation}) lies in
$K[\mathfrak{m}]$. But $A_m^{j_v,j_{w}}$ is contained in the ring
$K[\mathfrak{m}]$ provided
\raisebox{0pt}[13pt]{$A_m^{j_v,j_{w-1}}\in K[\mathfrak{m}]$}.
Finally from Proposition~\ref{C=C_as-minor} the function
$$A_m^{j_v,j_{v+1}}=C_{m}^{j_v}\cdot M_{\xi'}$$
lies in the ring $K[\mathfrak{m}]$.~$\Box$

\thenv{Conjecture\label{Hyp-small}}{Suppose the reductive
subalgebra~$\mathfrak{r}$ consists of three blocks; then the algebra
of invariants $K[\mathfrak{m}]^N$ is generated by the polynomials
$M_{\xi}$, $\xi\in S$, $L_{\varphi}$, $\varphi\in\Phi$, and the
elements (\ref{A}),~(\ref{B}),~(\ref{C}).}

%=======================================================================
%=======================================================================
%=======================================================================

\sect{The algebra of invariants in the case $(2,k,2)$}

Assume that the reductive subalgebra of the parabolic subalgebra
consists of three blocks with sizes $(2,k,2)$, $k>3$. We shall prove
Conjecture~\ref{Hyp-small} in this case.

The roots
$$S=\Big\{(1,4),(2,3),(k+1,k+4),(k+2,k+3)\Big\},$$
$$\Phi=\Big\{(3,k+3),(3,k+4),(4,k+3),(4,k+4)\Big\}$$
form the base and the system~$\Phi$, respectively. We have the
following diagram for the parabolic subalgebra in the case
$(2,k,2)$.
\begin{center}\refstepcounter{diagram}\label{Diag_2k2}
{\begin{tabular}{|p{0.1cm}|p{0.1cm}|p{0.1cm}|p{0.1cm}|p{0.1cm}|
p{0.1cm}|p{0.1cm}|p{0.1cm}|p{0.1cm}|p{0.1cm}|p{0.1cm}|l}
\multicolumn{2}{l}{{\small 1\hspace{5pt} 2
}}&\multicolumn{2}{l}{{\small 3\hspace{5pt}
4}}&\multicolumn{8}{r}{{\small $k+4\qquad$}}\\
\cline{1-11} \multicolumn{2}{|l|}{1}&&$\otimes$&&\multicolumn{2}{|l|}{$\ldots$}&&&&&{\small 1}\\
\cline{3-11} \multicolumn{2}{|r|}{1}&$\otimes$&&&\multicolumn{2}{|l|}{$\ldots$}&&&&&{\small 2}\\
\cline{1-11} \multicolumn{2}{|c|}{}&\multicolumn{7}{|l|}{1}&
$\times$&$\times$&{\small 3}\\
\cline{10-11}
\multicolumn{2}{|c|}{}&\multicolumn{7}{|l|}{\hspace{15pt}1}&$\times$&$\times$&{\small 4}\\
\cline{10-11}
\multicolumn{2}{|c|}{}&\multicolumn{7}{|l|}{\hspace{30pt}1}&&&\\
\cline{10-11} \multicolumn{2}{|c|}{}&\multicolumn{7}{|c|}{$\ldots$}&\multicolumn{2}{|l|}{$\ldots$}&\\
\cline{10-11}
\multicolumn{2}{|c|}{}&\multicolumn{7}{|r|}{1\hspace{30pt}\,}&&&\\
\cline{10-11}
\multicolumn{2}{|c|}{}&\multicolumn{7}{|r|}{1\hspace{15pt}\,}&&$\otimes$&{\small $k+1$}\\
\cline{10-11} \multicolumn{2}{|c|}{}&\multicolumn{7}{|r|}{1}&$\otimes$&&{\small$k+2$}\\
\cline{3-11} \multicolumn{9}{|c|}{}&\multicolumn{2}{|l|}{1}&{\small$k+3$}\\
\multicolumn{9}{|c|}{}&\multicolumn{2}{|r|}{1}&{\small$k+4$}\\
\cline{1-11} \multicolumn{11}{c}{Diagram \arabic{diagram}}\\
\end{tabular}}
\end{center}

Denote
$$M_1=M_{(2,3)},\ M_2=M_{(1,4)},\ N_1=M_{(k+2,k+3)},\ N_2=M_{(k+1,k+4)};$$
$$L_{11}=L_{\bigl((2,3),(k+2,k+3)\bigr)},\ L_{12}=L_{\bigl((2,3),(k+1,k+4)\bigr)},$$
$$L_{21}=L_{\bigl((1,4),(k+2,k+3)\bigr)},\ L_{22}=L_{\bigl((1,4),(k+1,k+4)\bigr)}.$$
The nonzero polynomials~(\ref{A}),~(\ref{B}) are not defined; the
invariant~(\ref{C}) is the following one:
\begin{equation}
D=C_{3}^{k+3}=\frac{L_{12}L_{21}-L_{11}L_{22}}{M_1N_1}.\label{2k2_D=C}
\end{equation}

We show that the invariants $M_i$, $N_j$, and $L_{ij}$, $i,j=1,2$,
do not generate the algebra of invariants $K[\mathfrak{m}]^N$.
Denote
$$\mathcal{B}_0=K[M_1,M_2,N_1,N_2,L_{11},L_{12},L_{21},L_{22}].$$
Assume the converse. Then $K[\mathfrak{m}]^N=\mathcal{B}_0$.
Therefore from Theorem~\ref{Th_local_field} it follows that there
exists a polynomial $f$ in 8 variables and there are numbers~$l_i$,
$i=1,\ldots,4$, such that
$$D=\frac{f(M_1,M_2,N_1,N_2,L_{11},L_{12},L_{21},L_{22})}%
{M_1^{l_1}M_2^{l_2}N_1^{l_3}N_2^{l_4}}.$$ Using~(\ref{2k2_D=C}), we
get the identity, which contradicts the algebraic independence of
the polynomials $M_i$, $N_j$, and $L_{ij}$ (see
Theorem~\ref{M-L_independ}).

\thenv{Theorem\label{Th_2k2}}{\emph{For any parabolic subalgebra
$\mathfrak{p}$ such that the reductive subalgebra of $\mathfrak{p}$
consists of three block with sizes} $(2,k,2)$, $k>3$, \emph{it
follows that the algebra of invariants $K[\mathfrak{m}]^N$ is
generated by the polynomials} $M_i$, $N_i$, $L_{ij}$, \emph{and}
$D$, \emph{where} $i,j=1,2$.}

\textsc{Proof.} Denote the set
$$\mathcal{B}=K[M_i,N_j,L_{ij},D]_{i,j=1,2}\subset K[\mathfrak{m}]^N.$$
We show that $K[\mathfrak{m}]^N\subset \mathcal{B}$.

As above, let $\mathcal{S}$ be the set of denominators generated by
the minors $M_i$, $N_i$, where $i=1,2$. From
Theorem~\ref{Th_local_field} it follows that the localization
$K[\mathfrak{m}]^N_\mathcal{S}$ coincides with the algebra of
Laurent polynomials
$$K\left[(M_i)^{\pm1},\
(N_i)^{\pm1},\ L_{ij}\right]_{i,j=1,2}.$$ If $f\in
K[\mathfrak{m}]^N$, then there exist natural numbers
$k_1,k_2,l_1,l_2$ such that
$$f\cdot M_1^{k_1}M_2^{k_2}N_1^{l_1}N_2^{l_2}\in \mathcal{B}_0.$$
Suppose for any minor~$M$ from the collection
$\left\{M_i,N_i\right\}_{i=1,2}$ we have that
$$\mbox{if }f\in K[\mathfrak{m}]\mbox{ and }M\cdot f\in
\mathcal{B},\mbox{ then }f\in\mathcal{B};$$ then $f\in\mathcal{B}$.

Let us prove that if $M\cdot f\in\mathcal{B}$, then
$f\in\mathcal{B}$ in the case $M=M_1$. The proof for the other cases
is similarly. Suppose $f\in K[\mathfrak{m}]^N$ and
$M_1f\in\mathcal{B}$. The polynomial $M_1f$ is equal to zero on
$\mathrm{Ann}\ M_1$. Let $a_i$, $b_i$, $c_{ij}$, $i,j=1,2$, be any
numbers in $K$. We construct a matrix
$$Q_1=\left(\begin{array}{ccccccccc}
0&0&a_2&0&0&\ldots&0&0&0\\
0&0&0&a_1&0&\ldots&0&0&0\\
0&0&0&0&0&\ldots&0&c_{11}&c_{12}\\
0&0&0&0&0&\ldots&0&c_{21}&c_{22}\\
0&0&0&0&0&\ldots&0&0&0\\
\ldots&\ldots&\ldots&\ldots&\ldots&\ldots&\ldots&\ldots&\ldots\\
0&0&0&0&0&\ldots&0&0&b_2\\
0&0&0&0&0&\ldots&0&b_1&0\\
0&0&0&0&0&\ldots&0&0&0\\
0&0&0&0&0&\ldots&0&0&0\\
\end{array}\right).$$
By direct calculation, we have
$$M_1(Q_1)=0,\ M_2(Q_1)=a_1a_2,\ N_1(Q_1)=b_1,\ N_2(Q_1)=-b_1b_2;$$
$$L_{11}(Q_1)=a_1c_{21},\ L_{12}(Q_1)=-a_1b_1c_{22},$$
$$L_{21}(Q_1)=a_1a_2c_{21},\ L_{22}(Q_1)=a_1a_2b_1c_{22},$$
$$D(Q_1)=a_1a_2(c_{12}c_{21}-c_{11}c_{22}).$$

Consider an algebra of polynomials
$$\mathcal{A}=K[u_1,u_2,v_1,v_2,w_{11},w_{12},w_{21},w_{22},z]$$
in 9 variables. Since the polynomial $fM_1$ is contained
in~$\mathcal{B}$, we have there exists a polynomial
$p(u_1,u_2,v_1,v_2,w_{11},w_{12},w_{21},w_{22},z)$ in the algebra
$\mathcal{A}$ such that
\begin{equation}
fM_1=p(M_1,M_2,N_1,N_2,L_{11},L_{12},L_{21},L_{22},D).\label{fM_1}
\end{equation}
Since $fM_1=0$ on $\mathrm{Ann}\,M_1$, we obtain $(fM_1)(Q_1)=0$. It
now follows that
$$\begin{array}{l}0=(fM_1)(Q_1)=p\Big(0,a_1a_2,b_1,-b_1b_2,\\
a_2c_{21},-a_2b_1c_{22},a_1a_2c_{21},
-a_1a_2b_1c_{22},a_1a_2(c_{12}c_{21}-c_{11}c_{22})\!\Big).\\
\end{array}$$
By $\mathcal{Z}$ denote a subset in $K^9$ that consists of the
following collections
$$\begin{array}{c}
\Bigl(0,a_1a_2,b_1,-b_1b_2,a_2c_{21},-a_2b_1c_{22},a_1a_2c_{21},
\qquad\qquad\qquad\\
-a_1a_2b_1c_{22},a_1a_2(c_{12}c_{21}-c_{11}c_{22})\Big),\mbox{
where }a_i,b_j,c_{ij}\in K.\\
\end{array}$$
It can easily be checked that the polynomials $u_1$ and
$w_{12}w_{21}-w_{11}w_{22}$ are equal to zero on~$\mathcal{Z}$. We
show that the ideal $\mathcal{I}_\mathcal{Z}=
\Big\{\varphi\in\mathcal{A}:\varphi(\mathcal{Z})=0\Big\}$ is
generated by $u_1$ and $w_{12}w_{21}-w_{11}w_{22}$.

Let the ideal $\mathcal{I}$ be generated by the polynomials
$u_1,w_{12}w_{21}-w_{11}w_{22}$. Note that $\mathcal{I}$ is a
relevant prime ideal. Indeed, the polynomial
$w_{12}w_{21}-w_{11}w_{22}$ is a irreducible one. Therefore the
algebra
\begin{equation}
\begin{array}{c}
\mathcal{A}/\mathcal{I}=\mathcal{A}/\langle u_1,w_{12}w_{21}-w_{11}w_{22}\rangle=\\
=K[u_2,v_1,v_2,w_{11},w_{12},w_{21},w_{22},z]/
\langle w_{12}w_{21}-w_{11}w_{22}\rangle\\
\end{array}\label{A/<>}
\end{equation}
is a domain of integrity. Hence the ideal~$\mathcal{I}$ is a
relevant prime one.

Since the polynomials $u_1,w_{12}w_{21}-w_{11}w_{22}$ vanish on the
set~$\mathcal{Z}$, it follows that
$\mathcal{I}\subset\mathcal{I}_{\mathcal{Z}}$ and
$\mathrm{Ann}\,\mathcal{I}\supset\mathcal{Z}$. The dimension of the
variety $\mathrm{Ann}\,\mathcal{I}$ is equal to the transcendence
degree of the quotient field of the algebra~(\ref{A/<>}) over~$K$,
i.e., $\dim\mathrm{Ann}\,\mathcal{I}=7$.

Evidently, $\dim\mathcal{Z}=7$, then
$\dim\mathrm{Ann}\,\mathcal{I}=\dim\mathcal{Z}$. Since
\mbox{$\mathrm{Ann}\,\mathcal{I}\supset\mathcal{Z}$}, we have
$\mathrm{Ann}\,\mathcal{I}=\overline{\mathcal{Z}}$. We show that
$\mathcal{I}=\mathcal{I}_{\mathcal{Z}}$. Suppose
$g\in\mathcal{I}_{\mathcal{Z}}$. By the Hilbert's theorem on zeros,
there exists a natural number~$N$ such that $g^N\in\mathcal{I}$.
Since $\mathcal{I}$ is a relevant prime ideal, we have
$g\in\mathcal{I}$. Hence,
$\mathcal{I}_{\mathcal{Z}}\subset\mathcal{I}$. Using
$\mathcal{I}\subset\mathcal{I}_{\mathcal{Z}}$, we obtain
$\mathcal{I}_{\mathcal{Z}}=\mathcal{I}$.

Thus since $p|_{\mathcal{Z}}=0$, it follows that $p\in\mathcal{I}$.
Therefore there exist polynomials $p_1$ and $p_2$ in the
algebra~$\mathcal{A}$ such that
$$p=p_1u_1+p_2(w_{12}w_{21}-w_{11}w_{22}).$$
Combining this with (\ref{fM_1}), we obtain
$$fM_1=p_1M_1+p_2\Big(L_{12}L_{21}-L_{11}L_{22}\Big)=p_1M_1+p_2M_1N_1D.$$
Further, we have $f=p_1+p_2N_1D$. So, $f\in\mathcal{B}$.~$\Box$

\thenv{Corrolary\label{Cor_2k2}}{\emph{The algebra of invariants
$K[\mathfrak{m}]^N$ is isomorphic onto the factor algebra
$$K[X_1,X_2,X_3,X_4,Y_1,Y_2,Y_3,Y_4,Z]/(X_2X_4Z-Y_2Y_3+Y_1Y_4)$$
in the case} $(2,k,2)$, $k>3$.}

%=============================================================
%=============================================================
%=============================================================
\sect{The algebra of invariants in the case $(1,2,2,1)$}

Suppose the reductive subalgebra~$\mathfrak{r}$ of~$\mathfrak{p}$
consists of the blocks with sizes $(1,2,2,1)$. We describe the
algebra of invariants in this case. The expanded base consists of
the following roots:
$$S=\{(1,2),(3,4),(2,5),(5,6)\},\quad\Phi=\{(2,4),(4,6)\}.$$
The diagram has the form
\begin{center}\refstepcounter{diagram}
{\begin{tabular}{|p{0.1cm}|p{0.1cm}|p{0.1cm}|p{0.1cm}|p{0.1cm}|
p{0.1cm}|c} \multicolumn{2}{l}{{\small 1\hspace{5pt}
2}}&\multicolumn{2}{l}{{\small 3\hspace{5pt}
4}}&\multicolumn{2}{l}{{\small 5\hspace{5pt} 6}}\\
\cline{1-6} 1&$\otimes$&&&&&{\small 1}\\
\cline{1-6} &\multicolumn{2}{|l|}{1}&$\times$&$\otimes$&&{\small 2}\\
\cline{4-6} &\multicolumn{2}{|r|}{1}&$\otimes$&&&{\small 3}\\
\cline{2-6} \multicolumn{3}{|c|}{}&\multicolumn{2}{|l|}{1}&$\times$&{\small 4}\\
\cline{6-6} \multicolumn{3}{|c|}{}&\multicolumn{2}{|r|}{1}&$\otimes$&{\small 5}\\
\cline{4-6} \multicolumn{5}{|c|}{}&1&{\small 6}\\
\cline{1-6} \multicolumn{6}{c}{Diagram \arabic{diagram}}\\
\end{tabular}}
\end{center}

The roots of the expanded base correspond to the following
invariants.
\begin{equation}
\begin{array}{c}
M_1=M_{(1,2)}=x_{1,2},\ M_2=M_{(3,4)}=x_{3,4},\\
M_3=M_{(2,5)}=\left|
\begin{array}{cc}
x_{2,4}&x_{2,5}\\
x_{3,4}&x_{3,5}
\end{array}\right|,\ M_4=M_{(5,6)}=x_{5,6};\\
L_1=L_{(2,4)}=x_{1,2}x_{2,4}+x_{1,3}x_{3,4},\\
L_2=L_{(4,6)}=x_{3,4}x_{4,6}+x_{3,5}x_{5,6}.
\end{array}\label{Invs1221}
\end{equation}

Consider a polynomial
$$D=x_{1,2}x_{2,4}x_{4,6}+x_{1,2}x_{2,5}x_{5,6}+
x_{1,3}x_{3,4}x_{4,6}+x_{1,3}x_{3,5}x_{5,6}.$$

\thenv{Lemma\label{L_1221}}{\emph{The polynomial $D$ is an
$N$-invariant}, $D$ \emph{is equal to the entry $(1,6)$ of the
matrix} $\mathbb{X}^3$, \emph{and we have}
\begin{equation}
M_2D=L_1L_2-M_1M_3M_4.\label{D_1221}
\end{equation}}

\textsc{Proof.} This lemma can be proved by direct
calculations.~$\Box$

\medskip
Note that the algebra of invariants $K[\mathfrak{m}]^N$ does not
coincide with the algebra
$$\mathcal{B}_0=K[M_1,M_2,M_3,M_4,
L_1,L_2].$$ Indeed, if $K[\mathfrak{m}]^N=\mathcal{B}_0$, then by
Theorem~\ref{Th_local_field}, there exists a polynomial
$f(t_1,\ldots,t_6)$ and there are natural numbers $l_i$,
$i=1,\ldots,4$, such that
$$D=\frac{f\Big(M_1,M_2,M_3,M_4,
L_1,L_2\Big)} {M_1^{\,l_1}M_2^{\,l_2}M_3^{\,l_3}M_4^{\,l_4}}.$$
Hence using the equation~(\ref{D_1221}), we obtain that the
invariants~(\ref{Invs1221}) are algebraically depended. This
contradicts Theorem~\ref{M-L_independ}.

\thenv{Theorem\label{Th1221}}{\emph{Let the reductive subalgebra
consists of the blocks} $(1,2,2,1)$. \emph{Then the algebra of
invariants is generated by the invariants}~(\ref{Invs1221})
\emph{and by the polynomial}~$D$.}

The proof is similarly to the proof of Theorem~\ref{Th_2k2}.

\thenv{Corollary\label{Cor_1221}}{\emph{The algebra of invariants
$K[\mathfrak{m}]^N$ is isomorphic onto the factor algebra
$$K[X_1,X_2,X_3,X_4,Y_1,Y_2,Z]/(X_2Z-Y_1Y_2+X_1X_3X_4)$$ in the
case} $(1,2,2,1)$.}

So, we have proved that the invariant $D$ is the minor of order 1 of
the matrix $\mathbb{X}^3$ in the case $(1,2,2,1)$ (see
Lemma~\ref{L_1221}). By Proposition~\ref{C=C_as-minor}, the
invariants $C^{\,i}_{m}$ can be determined as combined minors, which
are constructed by the formal matrices~$\mathbb{X}$ and
$\mathbb{X}^2$. Similarly, the basic invariants $M_{\xi}$, $\xi\in
S$, and $L_{\varphi}$, $\varphi\in\Phi$, are combined minors (see
Proposition~\ref{Def_L_phi}). Conjecture~\ref{conjecture-main}
generalizes all cases.

\medskip
Let $\mathfrak{p}$ be any parabolic subalgebra. Let
$I_1,I_2,\ldots,I_k$ and $J_1,J_2,\ldots,J_k$ be any collections of
rows and columns, respectively, such that we have the following
properties.
\begin{enumerate}
\item $|I_1|+|I_2|+\ldots+|I_k|=|J_1|+|J_2|+\ldots+|J_k|$.
\item Let $I$ be a set in the collection $\{I_1,I_2,\ldots,I_k,
J_1,J_2,\ldots,J_k\}$. If $\min I<i<\max I$, then $i\in I$.
\item For any $l=1,\ldots,k-1$ we have
$$\max I_l<\min I_{l+1},\ \min J_l>\max J_{l+1}.$$
\item Suppose for a positive root $(a,b)$ and a number $l$ such that
$1\leqslant l\leqslant k$, we have one of the following conditions:
\begin{itemize}
\item[a)] $a>\max I_l$ and $b\in J_{k-l+1}$,
\item[b)] $a\in I_l$ and $b<\min J_{k-l+1}$;
\end{itemize}
then the root $(a,b)$ is not contained in~$M$. In other words, the
minor $\mathbb{X}_{I_i}^{J_{k-i+1}}$ is bordered of zeros to the
right and to the below in the matrix~$\mathbb{X}$ for any
$i=1,\ldots,k$.
\end{enumerate}

We form a determinant, which is defined by the minors of the formal
matrix $\mathbb{X},\mathbb{X}^2,\ldots,\mathbb{X}^k$:
$$D_{I_1,I_2,\ldots,I_k}^{J_1,J_2,\ldots,J_k}=
\left|\begin{array}{cccc}
\mathbb{X}_{I_1}^{J_k}&(\mathbb{X}^2)_{I_1}^{J_{k-1}}&\ldots&(\mathbb{X}^k)_{I_1}^{J_1}\\
0&\mathbb{X}_{I_2}^{J_{k-1}}&\ldots&(\mathbb{X}^{k-1})_{I_2}^{J_1}\\
\ldots&\ldots&\ldots&\ldots\\
0&0&\ldots&\mathbb{X}_{I_k}^{J_1}\\
\end{array}\right|.$$

\thenv{Conjecture~\label{conjecture-main}}{The algebra of invariants
$K[\mathfrak{m}]^N$ is generated by the combined minors
$D_{I_1,I_2,\ldots,I_k}^{J_1,J_2,\ldots,J_k}$, where the conditions
1--4 hold for $I_1,I_2,\ldots,I_k$ and $J_1,J_2,\ldots,J_k$}

\textsc{Department of Mechanics and Mathematics, Samara State
University, Russia}\\ \emph{E-mail address}:
\verb"victoria.sevostyanova@gmail.com"

\end{document}